\documentclass{article}%
\usepackage{amsmath}
\usepackage{amsfonts}
\usepackage{amssymb}
\usepackage{graphicx}%
\providecommand{\U}[1]{\protect\rule{.1in}{.1in}}
\newtheorem{theorem}{Theorem}

\newtheorem{lemma}[theorem]{Lemma}

\begin{document}

\title{Blow-up result for a piezoelectric beams system with magnetic effects}
\author{Mohammad Kafini\\Department of Mathematics\\The Interdisciplinary Research Center in Construction\\and Building Materials\\King Fahd University of Petroleum and Minerals, KFUPM\\Dhahran 31261, Saudi Arabia\\E-mail: mkafini@kfupm.edu.sa}
\date{}
\maketitle

\begin{abstract}
The system under studying is for a piezoelectric beams system with magnetic
effects, frictional dampings and source terms. We use the concavity method to
study the competition of the dampings and the sources that leads to a blow-up
result for solutions with negative initial energy.

\textbf{Keywords} \textbf{and phrases: }Blow up, Negative initial energy,
Piezoelectric beams, Magnetic effects, Nonlinear sources.

\textbf{AMS Classification : }35B44, 35D30, 35L05, 35L15, 35L70,

\end{abstract}

\section{Introduction}

In this work, we are concerned with the following initial-boundary value
problem%
\begin{equation}%
\begin{tabular}
[c]{ll}%
$\rho\upsilon_{tt}-\alpha\upsilon_{xx}+\gamma\beta p_{xx}+\upsilon_{t}%
=f_{1}\left(  \upsilon,p\right)  ,$ & $(0,L)\times(0,T),$\\
$\mu p_{tt}-\beta p_{xx}+\gamma\beta\upsilon_{xx}+p_{t}=f_{2}\left(
\upsilon,p\right)  ,$ & $(0,L)\times(0,T),$\\
$\upsilon(x,0)=\upsilon_{0}(x),\quad\upsilon_{t}(x,0)=\upsilon_{1}(x),$ &
$x\in(0,L),$\\
$p(x,0)=p_{0}(x),\quad p_{t}(x,0)=p_{1}(x),$ & $x\in(0,L),$\\
$\upsilon(0,t)=\upsilon_{x}(L,t)=p(0,t)=p_{x}(L,t)=0,$ & $t\in(0,T),$%
\end{tabular}
\label{1}%
\end{equation}
where $\left(  x,t\right)  \in(0,L)\times(0,T)$ for $T>0.$ The function
$\upsilon$ represents the transverse displacement of the beam and $p$ is the
magnetic current or the total load of the electric displacement along the
transverse direction at each point $x$. The coefficients $\rho$, $\alpha$,
$\gamma,\mu$ and $\beta$ denote mass density per unit volume, elastic
stiffness, piezoelectric coefficient, magnetic permeability, and
impermeability, respectively. In addition, the following relationship%
\begin{equation}
\alpha_{1}=\alpha-\gamma^{2}\beta>0 \label{2}%
\end{equation}
was considered. The nonlinear sources $f_{1}$ and $f_{2}$ are functions
satisfying some conditions to be specified.

Piezoelectric materials have the valuable property of converting mechanical
energy to electro-magnetic energy, and vice versa. Piezoelectric actuators are
generally scalable, smaller, less expensive and more efficient than
traditional actuators, and hence, a competitive choice for many tasks,
particularly the control of flexible structures. They are used in civil,
industrial, automotive, aeronautic, and space structures. For more information
we refer to $\cite{1},\cite{2},\cite{3},\cite{4},$ for instance.

It is known in the literature that magnetic effects in the piezoelectric beams
is relatively small, and does not change the overall dynamics. Therefore, most
models for piezoelectric beams completely ignore the magnetic energy. These
models are known to be exponentially stabilizable by a mechanical feedback
controller in the energy space. For example, a beam of length $L$ and
thickness $h$, the system $\left(  \ref{3}\right)  $ is for undamped
Euler-Bernoulli beam model and the electrostatic or quasi-static assumptions
describe the stretching motion as%
\begin{equation}%
\begin{tabular}
[c]{l}%
$\rho\upsilon_{tt}-\alpha_{1}\upsilon_{xx}=0,\text{ \quad}$in$\text{
\ }(0,L)\times\mathbb{R}^{+}$\\
$\upsilon(x,0)=\upsilon_{0}(x),\quad\upsilon_{t}(x,0)=\upsilon_{1}(x),\quad
x\in(0,L),$\\
$\upsilon(0,t)=0,$ $\alpha_{1}\upsilon_{x}(L,t)=-\frac{\gamma V\left(
t\right)  }{h},\quad t>0.$%
\end{tabular}
\label{3}%
\end{equation}
Here, $V(t)$ denotes the voltage applied at the electrodes, and $\upsilon$
denotes the longitudinal displacement of the beam. The control system $\left(
\ref{3}\right)  $ is a well-posed boundary control system on a Hilbert space,
with norm corresponding to the system energy. It has been shown that it is
exactly controllable and hence exponentially stabilizable. In fact, by
choosing the feedback $V\left(  t\right)  =-\upsilon_{t}(L,t)$, the solutions
of the system become exponentially stable $\cite{5}$, for example.

On the other hand, including fully dynamic magnetic effects in a model for a
piezoelectric beam with voltage-driven electrodes. Morris and \"{O}zer
$\cite{6}$ used a variational approach to derive the differential equations
and boundary conditions that model a single piezoelectric beam with magnetic
effects given by%
\begin{equation}%
\begin{tabular}
[c]{l}%
$\rho\upsilon_{tt}-\alpha_{1}\upsilon_{xx}+\gamma\beta p_{xx}=0,$ $\ \ \text{
\ }(0,L)\times\mathbb{R}^{+}$\\
$\mu p_{tt}-\beta p_{xx}+\gamma\beta\upsilon_{xx}=0,$ $\ \ \text{
\ }(0,L)\times$ $\mathbb{R}^{+}$\\
$\upsilon(x,0)=\upsilon_{0}(x),\quad\upsilon_{t}(x,0)=\upsilon_{1}(x),\quad
x\in(0,L),$\\
$p(x,0)=p_{0}(x),\quad p_{t}(x,0)=p_{1}(x),\quad x\in(0,L),$\\
$\upsilon(0,t)=p(0,t)=\alpha\upsilon_{x}(L,t)-\gamma\beta p_{x}(L,t)=0,\quad
t>0.$\\
$\beta p_{x}(L,t)-\gamma\beta p_{x}(L,t)=-\frac{V\left(  t\right)  }{h},\quad
t>0.$%
\end{tabular}
\label{4}%
\end{equation}

It is also worth to mention the work of Ramos et al. $\cite{7}$ in which the
authors consider the system with a delay term acting in the damping of the
displacement. Using an energy-based approach, they proved the existence and
uniqueness of the solutions, the exponential stability of the solutions. For
damping acting on the magnetic current $p_{t}$ with nonlinear source terms
$f_{1}\left(  \upsilon,p\right)  $,$f_{2}\left(  \upsilon,p\right)  $ and
external nonlinear sources $h_{1}$ and $h_{2},$ Freitas et al. $\cite{8}$
studied the system of the form%
\begin{equation}%
\begin{tabular}
[c]{l}%
$\rho\upsilon_{tt}-\alpha_{1}\upsilon_{xx}+\gamma\beta p_{xx}+f_{1}\left(
\upsilon,p\right)  +\upsilon_{t}=h_{1},$ $\ \ \text{ \ }(0,L)\times T$\\
$\mu p_{tt}-\beta p_{xx}+\gamma\beta\upsilon_{xx}+f_{2}\left(  \upsilon
,p\right)  +\mu_{1}p_{t}(x,t)+\mu_{2}p_{t}(x,t-\tau)=h_{2},$ $\ \ \text{
\ }(0,L)\times$ $T$%
\end{tabular}
\label{5}%
\end{equation}
They proved that the system $\left(  \ref{5}\right)  $ is gradient and
asymptotically smooth, which as a consequence, implies the existence of a
global attractor, which is characterized as unstable manifold of the set of
stationary solutions. For more results to piezoelectric beam with magnetic
effects system in many different directions, we refer to $\cite{9}-\cite{13}.
$

For works that discussed the local existence and lower bound of blow-up time
to problems of a coupled nonlinear wave equations, we mention here the works
$\cite{14},\cite{15},\cite{16},\cite{17},$ e.g.

Our aim is to obtain blow up results with upper bound and lower bound of the
blow up time. The system under studying is for a piezoelectric beams system
with magnetic effects, frictional damping and source terms . The method used
is the concavity method introduced by Levine $\cite{18},\cite{19}$. In the
next section we present necessary conditions that needed to achieve our
result. In section 3, we state and prove our main result, the upper bound of
the blow up time. While in sections 4 we estimate the lower bound of that time.

\section{Preliminaries}

In this section we prepare some material needed in the proof of our main
result. For this reason we assume that

\noindent\textbf{(G1)} There exists a function $I(\upsilon,p)\geq0$ such that%

\[
\frac{\partial I}{\partial\upsilon}=f_{1}(\upsilon,p),\text{ \ \ \ \ }%
\frac{\partial I}{\partial p}=f_{2}(\upsilon,p).
\]
\textbf{(G2)} There exists a constant $\eta>2$ such that%

\[
\int_{0}^{L}[\upsilon f_{1}(\upsilon,p)+pf_{2}(\upsilon,p)-\eta I(\upsilon
,p)]dx\geq0.\quad\forall\upsilon,p\in H_{0}^{1}(0,L)
\]

\noindent\textbf{(G3) }There exists a constant $d>0$ such that%
\begin{align*}
\left\vert f_{1}(\chi,\phi)\right\vert  &  \leq d\left(  \left\vert
\chi\right\vert ^{\beta_{1}}+\left\vert \phi\right\vert ^{\beta_{2}}\right)
,\text{ \ \ \ \ }\forall(\chi,\phi)\in\mathbb{R}^{2},\\
\left\vert f_{2}(\chi,\phi)\right\vert  &  \leq d\left(  \left\vert
\chi\right\vert ^{\beta_{3}}+\left\vert \phi\right\vert ^{\beta_{4}}\right)
,\text{ \ \ \ \ }\forall(\chi,\phi)\in\mathbb{R}^{2},
\end{align*}
where%
\[
\beta_{i}\geq1,\text{ \ \ \ }(n-2)\beta_{i}\leq n,\text{ \ \ \ \ \ \ }%
i=1,2,3,4.
\]
Note that condition (G3) is necessary for the existence of a local solution to
$\left(  \ref{1}\right)  $. As an example of a function that satisfying the
other conditions (G1) and (G2) is%
\[
I(\upsilon,p):=\frac{a}{\eta}\left\vert \upsilon-p\right\vert ^{\eta},
\]
where we have%
\[
\frac{\partial I}{\partial\upsilon}=f_{1}(\upsilon,p)=a\left\vert
\upsilon-p\right\vert ^{\eta-2}(\upsilon-p),\ \ \ \ \ \ \frac{\partial
I}{\partial p}=f_{2}(\upsilon,p)=-a\left\vert \upsilon-p\right\vert ^{\eta
-2}(\upsilon-p).
\]
We define the Energy space to be%
\[
H=\left\{  H^{2}\left(  0,L\right)  \cap H_{\ast}^{1}\left(  0,L\right)
\times H_{\ast}^{1}\left(  0,L\right)  \times H^{2}\left(  0,L\right)  \cap
H_{\ast}^{1}\left(  0,L\right)  \times H_{\ast}^{1}\left(  0,L\right)
\right\}  ,
\]
where%
\[
H_{\ast}^{1}\left(  0,L\right)  =\left\{  g\in H^{1}\left(  0,L\right)
;\text{ }g\left(  0\right)  =g_{x}\left(  L\right)  =0\right\}  .
\]
The existence theorem can be stated as

\begin{theorem}
\emph{For any }$\left(  \upsilon_{0},\upsilon_{1},p_{0},p_{1}\right)  \in
H,$\emph{\ there is a unique strong solution of the system }$\left(
\ref{1}\right)  $\emph{.}
\end{theorem}

\noindent The proof can be established using the semigroup theory by combining
the results in $\cite{6},$ $\cite{14}$ and $\cite{20}.$

To make the system $\left(  \ref{1}\right)  $ more convenient, we will divide
the first equation by $\rho$ and the second one by $\mu$ but we keep the other
coefficients with the same name. The new system $\left(  \ref{1}\right)  $
becomes like%
\begin{equation}
\left\{
\begin{tabular}
[c]{ll}%
$\upsilon_{tt}-\alpha\upsilon_{xx}+\gamma\beta p_{xx}+\lambda_{1}\upsilon
_{t}=f_{1}\left(  \upsilon,p\right)  ,$ & $(0,L)\times(0,T),$\\
$p_{tt}-\beta p_{xx}+\gamma\beta\upsilon_{xx}+\lambda_{2}p_{t}=f_{2}\left(
\upsilon,p\right)  ,$ & $(0,L)\times(0,T),$%
\end{tabular}
\right.  \label{P}%
\end{equation}
together with the same initial and boundary conditions.

\noindent The energy functional associated to the above system is given by,
for $t\geq0,$
\begin{equation}
E(t):=\frac{1}{2}\left(  \alpha_{1}||\upsilon_{x}||_{2}^{2}+\beta
||\gamma\upsilon_{x}-p_{x}||_{2}^{2}+||\upsilon_{t}||_{2}^{2}+||p_{t}%
||_{2}^{2}\right)  -\int_{0}^{L}I(\upsilon,p)dx, \label{6}%
\end{equation}
where $\alpha_{1}$ defined in $\left(  \ref{2}\right)  $.

Consequently, we present the following lemma.

\begin{lemma}
\noindent\emph{The solution of }$\left(  \ref{P}\right)  $\emph{\ satisfies}
\begin{equation}
E^{\prime}\left(  t\right)  =-\lambda_{1}||\upsilon_{t}||_{2}^{2}-\lambda
_{2}||p_{t}||_{2}^{2}\leq0. \label{7}%
\end{equation}

\end{lemma}

\noindent\textbf{Proof. }Multiplying equation $\left(  \ref{P}\right)  _{1}$
by $\upsilon_{t},$ equation $\left(  \ref{P}\right)  _{2}$ by $p_{t}$ and
integrating by parts considering the boundary conditions over $(0,L)$ then
summing up we get%
\begin{align*}
&  \frac{1}{2}\frac{d}{dt}\left[  \left(  \alpha-\beta\gamma^{2}\right)
||\upsilon_{x}||_{2}^{2}+\beta||\gamma\upsilon_{x}-p_{x}||_{2}^{2}%
+||\upsilon_{t}||_{2}^{2}+||p_{t}||_{2}^{2}\right]  -\int_{0}^{L}\left(
\upsilon_{t}f_{1}+p_{t}f_{2}\right)  dx\\
&  =-\lambda_{1}||\upsilon_{t}||_{2}^{2}-\lambda_{2}||p_{t}||_{2}^{2}.
\end{align*}
From (G1), we notice that%
\[
\int_{0}^{L}\upsilon_{t}f_{1}dx+\int_{0}^{L}p_{t}f_{2}dx=\frac{d}{dt}\int%
_{0}^{L}I(\upsilon,p)dx.
\]
As $\alpha-\beta\gamma^{2}=\alpha_{1}>0$, the result follows.

Finally, we define%

\begin{equation}
F(t):=\frac{1}{2}\int_{0}^{L}\left(  \left\vert \upsilon(x,t)\right\vert
^{2}+\left\vert p(x,t)\right\vert ^{2}\right)  dx+\frac{1}{2}L\left(
t\right)  +\frac{1}{2}b(t+t_{0})^{2}, \label{8}%
\end{equation}
for $t_{0}\smallskip>0$ and $b>0$ to be chosen later and
\begin{align*}
L\left(  t\right)   &  :=\int_{0}^{t}\int_{0}^{L}\left(  \lambda_{1}\left\vert
\upsilon(x,s)\right\vert ^{2}+\lambda_{2}\left\vert p(x,s)\right\vert
^{2}\right)  dxds\\
&  +\left(  T-t\right)  \int_{0}^{L}\left(  \lambda_{1}\left\vert \upsilon
_{0}(x)\right\vert ^{2}+\lambda_{2}\left\vert p_{0}(x)\right\vert ^{2}\right)
dx.
\end{align*}
For the functional $L\left(  t\right)  $, we have the following lemma.

\begin{lemma}
\emph{Along the solution of }$\left(  \ref{P}\right)  $\emph{, the functional
}$L\left(  t\right)  $\emph{\ satisfies}%
\[
L^{\prime}\left(  t\right)  :=2\int_{0}^{t}\int_{0}^{L}\left(  \lambda
_{1}\upsilon\upsilon_{t}(x,s)+\lambda_{2}pp_{t}(x,s)\right)  dxds.
\]
$\emph{and}$%
\[
L^{\prime\prime}\left(  t\right)  :=2\int_{0}^{L}\left(  \lambda_{1}%
\upsilon\upsilon_{t}(x,t)+\lambda_{2}pp_{t}(x,t)\right)  dx.
\]

\end{lemma}

\noindent\textbf{Proof}. By a direct differentiation, we obtain%
\begin{align*}
L^{\prime}\left(  t\right)   &  =\int_{0}^{L}\left(  \lambda_{1}\left\vert
\upsilon(x,t)\right\vert ^{2}+\lambda_{2}\left\vert p(x,t)\right\vert
^{2}\right)  dx-\int_{0}^{L}\left(  \lambda_{1}\left\vert \upsilon
_{0}(x)\right\vert ^{2}+\lambda_{2}\left\vert p_{0}(x)\right\vert ^{2}\right)
dx\\
&  =\int_{0}^{L}\int_{0}^{t}\frac{\partial}{\partial s}\left[  \left(
\lambda_{1}\left\vert \upsilon(x,s)\right\vert ^{2}+\lambda_{2}\left\vert
p(x,s)\right\vert ^{2}\right)  \right]  dx\\
&  =2\int_{0}^{L}\int_{0}^{t}\left(  \lambda_{1}\upsilon\upsilon
_{t}(x,s)+\lambda_{2}pp_{t}(x,s)\right)  dsdx.
\end{align*}
The derivative of $L\left(  t\right)  $ again gives%
\begin{align*}
L^{\prime\prime}\left(  t\right)   &  =\frac{d}{dt}\left[  \int_{0}^{L}\left(
\lambda_{1}\left\vert \upsilon(x,t)\right\vert ^{2}+\lambda_{2}\left\vert
p(x,t)\right\vert ^{2}\right)  dx-\int_{0}^{L}\left(  \lambda_{1}\left\vert
\upsilon_{0}(x)\right\vert ^{2}+\lambda_{2}\left\vert p_{0}(x)\right\vert
^{2}\right)  dx\right] \\
&  =2\int_{0}^{L}\left(  \lambda_{1}\upsilon\upsilon_{t}(x,t)+\lambda
_{2}pp_{t}(x,t)\right)  dx.
\end{align*}

\section{ Blow up}

In this section we state and prove our main result.\newline

\begin{theorem}
\emph{Assume that (G1)--(G3) are hold. Then for any initial data }%
\[
(\upsilon_{0},p_{0},\upsilon_{1},p_{1})\in\left(  H^{2}\left(  0,L\right)
\cap H_{\ast}^{1}\left(  0,L\right)  \right)  ^{2}\times\left(  H_{\ast}%
^{1}\left(  0,L\right)  \right)  ^{2},
\]
\emph{\ \ satisfying}%
\[
E(0)=\frac{1}{2}\left(  \alpha_{1}||\upsilon_{0x}||_{2}^{2}+\beta
||\gamma\upsilon_{0x}-p_{0x}||_{2}^{2}+||\upsilon_{1}||_{2}^{2}+||p_{1}%
||_{2}^{2}\right)  -\int_{0}^{L}I(\upsilon_{0},p_{0})dx<0,
\]
\emph{the corresponding solution of }$\left(  \ref{P}\right)  $\emph{\ blows
up in a finite time.}
\end{theorem}

\noindent Before the proof, we need the following technical lemmas.

\begin{lemma}
\emph{Along the solution of }$\left(  \ref{P}\right)  $\emph{, we estimate}%
\begin{align}
\int_{0}^{L}\left(  \rho\upsilon\upsilon_{tt}+\mu pp_{tt}\right)  dx  &
\geq\alpha_{1}\left(  \frac{\eta}{2}-1\right)  \left\Vert \upsilon
_{x}\right\Vert _{2}^{2}+\beta\left(  \frac{\eta}{2}-1\right)  ||\gamma
\upsilon_{x}-p_{x}||_{2}^{2}\nonumber\\
&  -\int_{0}^{L}\left(  \lambda_{1}\upsilon\upsilon_{t}+\lambda_{2}%
pp_{t}\right)  dx+\frac{\eta}{2}\left(  ||\upsilon_{t}||_{2}^{2}+||p_{t}%
||_{2}^{2}\right) \nonumber\\
&  -\eta E(t). \label{9}%
\end{align}

\end{lemma}

\noindent\textbf{Proof. }Multiply the equations in $\left(  \ref{P}\right)
_{1}$ by $\upsilon,$ equation $\left(  \ref{P}\right)  _{2}$ by $p$ and
integrate by parts over $\left(  0,L\right)  $ to get%

\begin{align*}
\int_{0}^{L}\left(  \upsilon\upsilon_{tt}+pp_{tt}\right)  dx  &  =-\alpha
_{1}\left\Vert \upsilon_{x}\right\Vert ^{2}-\beta||\gamma\upsilon_{x}%
-p_{x}||_{2}^{2}-\int_{0}^{L}\left(  \lambda_{1}\upsilon\upsilon_{t}%
+\lambda_{2}pp_{t}\right)  dx\\
&  +\int_{0}^{L}[\upsilon f_{1}(\upsilon,p)+pf_{2}(\upsilon,p)]dx.
\end{align*}
Considering \textbf{(G2), }we have%
\begin{align*}
\int_{0}^{L}\left(  \upsilon\upsilon_{tt}+pp_{tt}\right)  dx  &  \geq
-\alpha_{1}||\upsilon_{x}||_{2}^{2}-\beta||\gamma\upsilon_{x}-p_{x}||_{2}%
^{2}\\
&  -\int_{0}^{L}\left(  \lambda_{1}\upsilon\upsilon_{t}+\lambda_{2}%
pp_{t}\right)  dx+\eta\int_{0}^{L}I(\upsilon,p)dx.
\end{align*}

\noindent Exploiting $\left(  \ref{6}\right)  $, we find%
\[
\eta\int_{0}^{L}I(\upsilon,p)dx=-\eta E(t)+\frac{\eta}{2}\left[  \alpha
_{1}||\upsilon_{x}||_{2}^{2}+\beta||\gamma\upsilon_{x}-p_{x}||_{2}%
^{2}+||\upsilon_{t}||_{2}^{2}+||p_{t}||_{2}^{2}\right]  ,
\]
to reach the result%
\begin{align*}
\int_{0}^{L}\left(  \upsilon\upsilon_{tt}+pp_{tt}\right)  dx  &  \geq
\alpha_{1}\left(  \frac{\eta}{2}-1\right)  ||\upsilon_{x}||_{2}^{2}%
+\beta\left(  \frac{\eta}{2}-1\right)  ||\gamma\upsilon_{x}-p_{x}||_{2}^{2}\\
&  -\int_{0}^{L}\left(  \lambda_{1}\upsilon\upsilon_{t}+\lambda_{2}%
pp_{t}\right)  dx+\frac{\eta}{2}\left(  ||\upsilon_{t}||_{2}^{2}+||p_{t}%
||_{2}^{2}\right)  -\eta E(t).
\end{align*}

\noindent This completes the proof.

\begin{lemma}
\noindent\textbf{\ }\emph{Along the solution of }$\left(  \ref{P}\right)
$\emph{, we estimate, for any }$\varepsilon>0,$%
\begin{align}
I  &  =\left(  \int_{0}^{L}\left(  \upsilon\upsilon_{t}+pp_{t}\right)
dx+\frac{1}{2}L^{\prime}\left(  t\right)  +b(t+t_{0})\right)  ^{2}\label{17}\\
&  \leq F(t)\left[  (1+\varepsilon)\left(  \int_{0}^{L}\upsilon_{t}^{2}%
dx+\int_{0}^{L}p_{t}^{2}dx\right)  +2(1+\frac{1}{\varepsilon})\left(
b-2\lambda E\left(  t\right)  \right)  \right]  .\nonumber
\end{align}

\end{lemma}

\noindent\textbf{Proof. }By using Young's inequality, we have%
\begin{align*}
I  &  =\left(  \int_{0}^{L}\left(  \upsilon\upsilon_{t}+pp_{t}\right)
dx\right)  ^{2}+2\left(  \frac{1}{2}L^{\prime}\left(  t\right)  +b(t+t_{0}%
)\right)  \int_{0}^{L}\left(  \upsilon\upsilon_{t}+pp_{t}\right)  dx\\
&  +\left(  \frac{1}{2}L^{\prime}\left(  t\right)  +b(t+t_{0})\right)  ^{2}\\
&  \leq\left(  \int_{0}^{L}\left(  \upsilon\upsilon_{t}+pp_{t}\right)
dx\right)  ^{2}+2\left[  \frac{\varepsilon}{2}\left(  \int_{0}^{L}\left(
\upsilon\upsilon_{t}+pp_{t}\right)  dx\right)  ^{2}+\frac{1}{2\varepsilon
}\left(  \frac{1}{2}L^{\prime}\left(  t\right)  +b(t+t_{0})\right)
^{2}\right] \\
&  +\left(  \frac{1}{2}L^{\prime}\left(  t\right)  +b(t+t_{0})\right)  ^{2}\\
&  \leq(1+\varepsilon)\left(  \int_{0}^{L}\left(  \upsilon\upsilon_{t}%
+pp_{t}\right)  dx\right)  ^{2}+(1+\frac{1}{\varepsilon})\left(  \frac{1}%
{2}L^{\prime}\left(  t\right)  +b(t+t_{0})\right)  ^{2}.
\end{align*}

Using Cauchy-Schwartz inequality yields%
\begin{align*}
I  &  \leq(1+\varepsilon)\left(  \int_{0}^{L}\upsilon^{2}dx+\int_{0}^{L}%
p^{2}dx\right)  \left(  \int_{0}^{L}\upsilon_{t}^{2}dx+\int_{0}^{L}p_{t}%
^{2}dx\right) \\
&  +2(1+\frac{1}{\varepsilon})\left[  \left(  \frac{1}{2}L^{\prime}\left(
t\right)  \right)  ^{2}+b^{2}(t+t_{0})^{2}\right]  .
\end{align*}
Using Lemma 2.2, we estimate%
\begin{align*}
&  \left(  \frac{1}{2}L^{\prime}\left(  t\right)  \right)  ^{2}\\
&  =\left(  \int_{0}^{t}\int_{0}^{L}\left(  \lambda_{1}\upsilon\upsilon
_{t}(x,s)+\lambda_{2}pp_{t}(x,s)\right)  dxds\right)  ^{2}\\
&  =\left(  \int_{0}^{t}\left(  \lambda_{1}||\upsilon\left(  s\right)
||_{2}^{2}+\lambda_{2}||p\left(  s\right)  ||_{2}^{2}\right)  ds\right)
\left(  \int_{0}^{t}\left(  \lambda_{1}||\upsilon_{t}\left(  s\right)
||_{2}^{2}+\lambda_{2}||p_{t}\left(  s\right)  ||_{2}^{2}\right)  ds\right)  .
\end{align*}
Using $\left(  \ref{7}\right)  ,$ we find%
\[
\int_{0}^{t}\left(  \lambda_{1}||\upsilon_{t}\left(  s\right)  ||_{2}%
^{2}+\lambda_{2}||p_{t}\left(  s\right)  ||_{2}^{2}\right)  ds\leq-E\left(
t\right)  ,
\]
together with $\left(  \ref{8}\right)  $ imply that
\[
\left(  \frac{1}{2}L^{\prime}\left(  t\right)  \right)  ^{2}\leq-2\lambda
F(t)E\left(  t\right)  .
\]
Hence%
\[
I\leq F(t)\left[  (1+\varepsilon)\left(  \int_{0}^{L}\upsilon_{t}^{2}%
dx+\int_{0}^{L}p_{t}^{2}dx\right)  +2(1+\frac{1}{\varepsilon})\left(
b-2\lambda E\left(  t\right)  \right)  \right]  .
\]
This completes the proof.

\noindent\noindent\textbf{Proof of theorem 3.1.} From $\left(  \ref{7}\right)
$, we notice that
\begin{equation}
E(t)\leq E(0)\leq0. \label{11}%
\end{equation}
By differentiating $F$ defined in $\left(  \ref{8}\right)  $ twice we get%

\begin{equation}
F^{\prime}(t)=\int_{0}^{L}\left(  \upsilon\upsilon_{t}+pp_{t}\right)
dx+\frac{1}{2}L^{\prime}\left(  t\right)  +b(t+t_{0}) \label{12}%
\end{equation}
and%

\begin{equation}
F^{\prime\prime}(t)=\int_{0}^{L}\left(  \upsilon\upsilon_{tt}+pp_{tt}\right)
dx+\int_{0}^{L}\left(  \left\vert \upsilon_{t}\right\vert ^{2}+\left\vert
p_{t}\right\vert ^{2}\right)  dx+\frac{1}{2}L^{\prime\prime}\left(  t\right)
+b. \label{13}%
\end{equation}
Using the estimation in $\left(  \ref{9}\right)  $, yields%
\begin{align}
F^{\prime\prime}(t)  &  \geq\alpha_{1}\left(  \frac{\eta}{2}-1\right)
\left\Vert \upsilon_{x}\right\Vert _{2}^{2}++\beta\left(  \frac{\eta}%
{2}-1\right)  ||\gamma\upsilon_{x}-p_{x}||_{2}^{2}\label{14}\\
&  +\left(  \frac{\eta}{2}+1\right)  \left\Vert \upsilon_{t}\right\Vert
_{2}^{2}+\left(  \frac{\eta}{2}+1\right)  \left\Vert p_{t}\right\Vert _{2}%
^{2}-\eta E\left(  t\right)  +b.\nonumber
\end{align}

We then define%
\begin{equation}
G(t):=F^{-\sigma}(t), \label{15}%
\end{equation}
for $\sigma>0$ to be chosen properly.

By differentiating $G$ twice we arrive at%
\begin{align*}
G^{\prime}(t)  &  =-\sigma F^{-(\sigma+1)}(t)F^{\prime}(t),\\
G^{\prime\prime}(t)  &  =-\sigma F^{-(\sigma+2)}(t)Q(t),
\end{align*}
where%
\begin{equation}
Q(t)=F(t)F^{\prime\prime}(t)-(\sigma+1)\left(  F^{\prime}(t)\right)  ^{2}.
\label{16}%
\end{equation}
Using $\left(  \ref{17}\right)  $, the last term in $\left(  \ref{16}\right)
$ becomes%
\begin{align}
&  -(\sigma+1)\left(  F^{\prime}(t)\right)  ^{2}\label{18}\\
&  =-(\sigma+1)\left(  \int_{0}^{L}(\upsilon\upsilon_{t}+pp_{t})dx+\frac{1}%
{2}L^{\prime}\left(  t\right)  +b(t+t_{0})\right)  ^{2}\nonumber\\
&  \geq-(\sigma+1)F(t)\left[  (1+\varepsilon)\left(  \int_{0}^{L}\upsilon
_{t}^{2}dx+\int_{0}^{L}p_{t}^{2}dx\right)  +2(1+\frac{1}{\varepsilon})\left(
b-2\lambda E\left(  t\right)  \right)  \right]  .\nonumber
\end{align}
Inserting $\left(  \ref{14}\right)  $ and $\left(  \ref{18}\right)  $ in
$\left(  \ref{16}\right)  $ yields%
\begin{align}
Q(t)  &  =F(t)F^{\prime\prime}(t)-(\sigma+1)\left(  F^{\prime}(t)\right)
^{2}\nonumber\\
&  \geq F(t)\alpha_{1}\left(  \frac{\eta}{2}-1\right)  \left\Vert \upsilon
_{x}\right\Vert _{2}^{2}+F(t)\beta\left(  \frac{\eta}{2}-1\right)
||\gamma\upsilon_{x}-p_{x}||_{2}^{2}\nonumber\\
&  +F(t)\left(  \frac{\eta}{2}+1-(\sigma+1)(1+\varepsilon)\right)  \left\Vert
\upsilon_{t}\right\Vert _{2}^{2}\label{181}\\
&  +F(t)\left(  \frac{\eta}{2}+1-(\sigma+1)(1+\varepsilon)\right)  \left\Vert
p_{t}\right\Vert _{2}^{2}\nonumber\\
&  +F(t)\left[  -\eta E\left(  t\right)  -2(\sigma+1)\left(  1+\frac
{1}{\varepsilon}\right)  \left(  b-2\lambda E\left(  t\right)  \right)
\right]  ,\text{\ \ \ \ \ \ \ }\forall\varepsilon>0.\nonumber
\end{align}

\noindent As $\eta$ $>2,$ we choose $\varepsilon$ satisfying%
\[
0<\varepsilon<\frac{\eta}{2}\text{ }%
\]
and%
\[
0<\sigma<\frac{\eta-2\varepsilon}{2\left(  1+\varepsilon\right)  }.
\]
Hence, $\left(  \ref{181}\right)  $ takes the form%
\[
Q(t)\geq F(t)\left\{  -E\left(  t\right)  \left[  \eta-4\lambda(\sigma
+1)\left(  1+\frac{1}{\varepsilon}\right)  \right]  -2b(\sigma+1)\left(
1+\frac{1}{\varepsilon}\right)  \right\}  ,\text{\ \ }\forall\varepsilon>0.
\]
At this point we may choose $\eta$ more large, if needed, that makes
\begin{equation}
\eta-4\lambda(\sigma+1)\left(  1+\frac{1}{\varepsilon}\right)  =k>0\text{.}
\label{182}%
\end{equation}
Finally, as $-kE\left(  0\right)  >0,$ we pick $b>0$ small enough so that%

\[
-kE\left(  0\right)  -2b(\sigma+1)\left(  1+\frac{1}{\varepsilon}\right)
\geq0.
\]

\noindent Therefore,
\[
Q(t)\geq0,\text{ \ \ \ }\forall t\geq0
\]

\noindent and hence
\[
G^{\prime\prime}(t)\leq0,\quad\forall t\geq0;
\]
that is $G^{\prime}$ is decreasing.

By choosing $t_{0}>0$ large enough we get%

\[
F^{\prime}(0)=\int_{0}^{L}(\upsilon_{0}\upsilon_{1}+p_{0}p_{1})dx+bt_{0}>0,
\]
hence $G^{\prime}(0)<0.$

Finally Taylor expansion of $G$ yields%

\[
G(t)\leq G(0)+tG^{\prime}(0),\quad\forall t\geq0,
\]
which shows that $G(t)$ vanishes at a time $t_{m}\leq-\frac{G(0)}{G^{\prime
}(0)}.$ Consequently $F(t)$ must become infinite at time $t_{m}.$

Therefore, the solution $(\upsilon,p)$ blows up with time bounded above by
$t_{m}.$

\section{\noindent Lower bound for the blow-up time}

\noindent In this section, we estimate a lower bound for the blow-up time. For
that reason, we define
\begin{equation}
\psi(t)=\int_{0}^{L}I(\upsilon,p)dx \label{183}%
\end{equation}
and have the following theorem.

\begin{theorem}
\emph{Assume that }$(\upsilon,p)$ \emph{is a solution of }$\left(
\ref{P}\right)  $\emph{\ that satisfying (G1)-(G3) and does not exist after a
time }$T^{\ast}$\emph{. Then we have}
\[
\int_{\psi(0)}^{\infty}\frac{c\text{ }dy}{y+y^{\beta_{1}}+y^{\beta_{2}%
}+y^{\beta_{3}}+y^{\beta_{4}}}\leq T^{\ast},
\]
\emph{for some constant }$c>0.$
\end{theorem}

\noindent\textbf{Proof.} A direct differentiation of $\psi(t)$ gives%
\[
\psi^{\prime}(t)=\int_{0}^{L}\left[  \upsilon_{t}I_{\upsilon}(\upsilon
,p)+p_{t}I_{p}(\upsilon,p)\right]  dx.
\]
Using Young's inequality and \textbf{(G1)}, we have%
\begin{align}
\psi^{\prime}(t)  &  \leq\frac{1}{2}\int_{0}^{L}\left(  \upsilon_{t}^{2}%
+p_{t}^{2}\right)  dx+\frac{1}{2}\int_{0}^{L}\left(  I_{\upsilon}^{2}%
+I_{p}^{2}\right)  dx\nonumber\\
&  =\frac{1}{2}\int_{0}^{L}\left(  \upsilon_{t}^{2}+p_{t}^{2}\right)
dx+\frac{1}{2}\int_{0}^{L}\left(  f_{1}^{2}+f_{2}^{2}\right)  dx. \label{19}%
\end{align}
Using \textbf{(G3) }gives%
\[
f_{1}^{2}(\upsilon,p)\leq d^{2}\left(  \left\vert \upsilon\right\vert
^{\beta_{1}}+\left\vert p\right\vert ^{\beta_{2}}\right)  ^{2}\leq
2d^{2}\left(  \left\vert \upsilon\right\vert ^{2\beta_{1}}+\left\vert
p\right\vert ^{2\beta_{2}}\right)
\]
\ \ \ and \ \
\[
\ \ f_{2}^{2}(\upsilon,p)\leq d^{2}\left(  \left\vert \upsilon\right\vert
^{\beta_{3}}+\left\vert p\right\vert ^{\beta_{4}}\right)  ^{2}\leq
2d^{2}\left(  \left\vert \upsilon\right\vert ^{2\beta_{3}}+\left\vert
p\right\vert ^{2\beta_{4}}\right)  .
\]
Using the embedding and Poincar\`{e}'s inequality, we obtain from $\left(
\ref{19}\right)  ,$
\begin{align*}
\psi^{\prime}(t)  &  \leq\frac{1}{2}\int_{0}^{L}\left(  \upsilon_{t}^{2}%
+p_{t}^{2}\right)  dx+d^{2}\int_{0}^{L}\left(  \left\vert \upsilon\right\vert
^{2\beta_{1}}+\left\vert p\right\vert ^{2\beta_{2}}+\left\vert \upsilon
\right\vert ^{2\beta_{3}}+\left\vert p\right\vert ^{2\beta_{4}}\right)  dx\\
&  \leq\frac{1}{2}\int_{0}^{L}\left(  \upsilon_{t}^{2}+p_{t}^{2}\right)  dx\\
&  +d^{2}\left(  C_{p}^{\beta_{1}}||\upsilon_{x}||_{2}^{2\beta_{1}}%
+C_{p}^{\beta_{2}}||p_{x}||_{2}^{2\beta_{2}}+C_{p}^{\beta_{3}}||\upsilon
_{x}||_{2}^{2\beta_{3}}+C_{p}^{\beta_{4}}||p_{x}||_{2}^{2\beta_{4}}\right)  ,
\end{align*}
where $C_{p}$ is the Poincar\`{e}'s constant.

Thus, for $a=d^{2}\max\left\{  C_{p}^{\beta_{1}},C_{p}^{\beta_{2}}%
,C_{p}^{\beta_{3}},C_{p}^{\beta_{4}}\right\}  ,$ we get%
\begin{equation}
\psi^{\prime}(t)\leq\frac{1}{2}\left(  \left\Vert \upsilon_{t}\right\Vert
_{2}^{2}+\left\Vert p_{t}\right\Vert _{2}^{2}\right)  +a\left(  ||\upsilon
_{x}||_{2}^{2\beta_{1}}+||p_{x}||_{2}^{2\beta_{2}}+||\upsilon_{x}%
||_{2}^{2\beta_{3}}+||p_{x}||_{2}^{2\beta_{4}}\right)  . \label{190}%
\end{equation}

\noindent Recalling $\left(  \ref{6}\right)  $ and $\left(  \ref{11}\right)  $
we infer that%
\begin{equation}
\alpha_{1}||\upsilon_{x}||_{2}^{2}+\beta||\gamma\upsilon_{x}-p_{x}||_{2}%
^{2}+\rho||\upsilon_{t}||_{2}^{2}+\mu||p_{t}||_{2}^{2}\leq2\int_{0}%
^{L}I(\upsilon,p)dx. \label{191}%
\end{equation}
If we let $m=\min\left\{  \alpha_{1},\beta,\rho,\mu\right\}  $ and use
$\left(  \ref{183}\right)  $ then $\left(  \ref{191}\right)  $ becomes%
\[
m\left(  ||\upsilon_{x}||_{2}^{2}+||\gamma\upsilon_{x}-p_{x}||_{2}%
^{2}+||\upsilon_{t}||_{2}^{2}+||p_{t}||_{2}^{2}\right)  \leq2\int_{0}%
^{L}I(\upsilon,p)dx=2\psi(t).
\]
So, each term in $\left(  \ref{190}\right)  $ is bounded above by $\frac{2}%
{m}\psi(t).$

Therefore, $\left(  \ref{190}\right)  $ yields
\begin{equation}
\psi^{\prime}(t)\leq\frac{2}{m}\psi(t)+a\left(  \frac{2}{m}\right)
^{r}\left[  \psi^{\beta_{1}}(t)+\psi^{\beta_{2}}(t)+\psi^{\beta_{3}}%
(t)+\psi^{\beta_{4}}(t)\right]  , \label{20}%
\end{equation}
where $r=\max\left\{  \beta_{1},\beta_{2},\beta_{3},\beta_{4}\right\}  .$

If we integrate $\left(  \ref{20}\right)  $ over $(0,T^{\ast})$ and use the
substitution $y=\psi(t),$ then we reach%
\[
\int_{\psi(0)}^{\infty}\frac{dy}{\frac{2}{m}y+a\left(  \frac{2}{m}\right)
^{r}\left[  y^{\beta_{1}}+y^{\beta_{2}}+y^{\beta_{3}}+y^{\beta_{4}}\right]
}\leq T^{\ast}.
\]
This completes the proof.

\begin{description}
\item[\textbf{Conclusion:}] In our proof of the blow-up property in theorem
3.1, we assume the domain $\left[  0,L\right]  $ is bounded. We notify here
that the same proof is also valid in unbounded domain. Moreover, condition
$\left(  \ref{182}\right)  $ is assumed to be a compitition between the source
(assumed $\eta$ large) and the damping (where it may be assumed small by
taking $\lambda$ small). Finally, positive initial energy can be considered
but with different arrangement.\noindent

\item[\textbf{Acknowledgment:}] The author would like to express his sincere
thanks to King Fahd University of Petroleum and Minerals and The
Interdisciplinary Research Center in Construction and Building Materials, for
support. This work has been funded by KFUPM under Project \# INCB2402.
\end{description}

\end{document}